\documentclass{wscpaperproc}
\usepackage{latexsym}
\usepackage{graphicx}
\usepackage{mathptmx}
\usepackage[T1]{fontenc}
\usepackage{float}
\usepackage{subfigure}
\usepackage{subcaption}
\usepackage{arydshln}
\usepackage{booktabs}
\usepackage{multirow} 
\usepackage{multicol} 
\usepackage{array}
\usepackage{indentfirst}
\usepackage{amsmath}
\usepackage{amsfonts}
\usepackage{amssymb}
\usepackage{amsbsy}
\usepackage{amsthm}
\usepackage[dvipsnames]{xcolor}
\usepackage{color}
\usepackage[markup=underlined, commentmarkup=todo, deletedmarkup=sout
, final
]{changes}
\usepackage[pdftex,colorlinks=true,urlcolor=blue,citecolor=black,anchorcolor=black,linkcolor=black]{hyperref}

\definechangesauthor[name={Yang}, color=RedOrange]{Yang} 
\definechangesauthor[name={Cao}, color=RoyalBlue]{Cao} 
\definechangesauthor[name={Cao1}, color=Purple]{Cao1} 
\definechangesauthor[name={Liang}, color=RedViolet]{Liang} 

\newtheoremstyle{wsc}
{3pt}
{3pt}
{}
{}
{\bf}
{}
{.5em}
{}

\theoremstyle{wsc}

\begin{document}

%
%

\pagestyle{fancyplain}

\thispagestyle{plain}
\firstPageHead{}

\chead{\fancyplain{}{\itshape Cao, Yang, Leong, Liu and Chan}}

\setlength{\headheight}{26pt}

\rhead{\thepage}
\cfoot{}
\renewcommand{\headrulewidth}{0pt} 

\makeatletter
\let\@internalcite\cite
\def\cite{\def\@citeseppen{-1000}%
    \def\@cite##1##2{(##1\if@tempswa , ##2\fi)}%
    \def\citeauthoryear##1##2##3{##1 ##3}\@internalcite}
\def\citeNP{\def\@citeseppen{-1000}%
    \def\@cite##1##2{##1\if@tempswa , ##2\fi}%
    \def\citeauthoryear##1##2##3{##1 ##3}\@internalcite}
\def\citeN{\def\@citeseppen{-1000}%
    \def\@cite##1##2{##1\if@tempswa, ##2)\else{}\fi}%
    \def\citeauthoryear##1##2##3{##1 (##3)}\@citedata}
\def\citeA{\def\@citeseppen{-1000}%
    \def\@cite##1##2{(##1\if@tempswa , ##2\fi)}%
    \def\citeauthoryear##1##2##3{##1}\@internalcite}
\def\citeANP{\def\@citeseppen{-1000}%
    \def\@cite##1##2{##1\if@tempswa , ##2\fi}%
    \def\citeauthoryear##1##2##3{##1}\@internalcite}
\def\shortcite{\def\@citeseppen{-1000}%
    \def\@cite##1##2{(##1\if@tempswa , ##2\fi)}%
    \def\citeauthoryear##1##2##3{##2 ##3}\@internalcite}
\def\shortciteNP{\def\@citeseppen{-1000}%
    \def\@cite##1##2{##1\if@tempswa , ##2\fi}%
    \def\citeauthoryear##1##2##3{##2 ##3}\@internalcite}
\def\shortciteN{\def\@citeseppen{-1000}%
    \def\@cite##1##2{##1\if@tempswa, ##2\else{}\fi}%
    \def\citeauthoryear##1##2##3{##2 (##3)}\@citedata}
\def\shortciteA{\def\@citeseppen{-1000}%
    \def\@cite##1##2{(##1\if@tempswa , ##2\fi)}%
    \def\citeauthoryear##1##2##3{##2}\@internalcite}
\def\shortciteANP{\def\@citeseppen{-1000}%
    \def\@cite##1##2{##1\if@tempswa , ##2\fi}%
    \def\citeauthoryear##1##2##3{##2}\@internalcite}
\def\citeyear{\def\@citeseppen{-1000}%
    \def\@cite##1##2{(##1\if@tempswa , ##2\fi)}%
    \def\citeauthoryear##1##2##3{##3}\@citedata}
\def\citeyearNP{\def\@citeseppen{-1000}%
    \def\@cite##1##2{##1\if@tempswa , ##2\fi}%
    \def\citeauthoryear##1##2##3{##3}\@citedata}
%
%
%
\def\@citedata{%
    \@ifnextchar [{\@tempswatrue\@citedatax}%
                  {\@tempswafalse\@citedatax[]}%
}

\def\@citedatax[#1]#2{%
\if@filesw\immediate\write\@auxout{\string\citation{#2}}\fi%
  \def\@citea{}\@cite{\@for\@citeb:=#2\do%
    {\@citea\def\@citea{, }\@ifundefined
       {b@\@citeb}{{\bf ?}%
       \@warning{Citation `\@citeb' on page \thepage \space undefined}}%
{\csname b@\@citeb\endcsname}}}{#1}}%

%
\def\@citex[#1]#2{%
\if@filesw\immediate\write\@auxout{\string\citation{#2}}\fi%
  \def\@citea{}\@cite{\@for\@citeb:=#2\do%
    {\@citea\def\@citea{; }\@ifundefined
       {b@\@citeb}{{\bf ?}%
       \@warning{Citation `\@citeb' on page \thepage \space undefined}}%
{\csname b@\@citeb\endcsname}}}{#1}}%

%
\def\@biblabel#1{}
\makeatother



\newdimen\bibindent
\bibindent=0.0em
\def\thebibliography#1{\section*{\refname}\list
   {}{\settowidth\labelwidth{[#1]}
   \leftmargin\parindent
   \itemindent -\parindent
   \listparindent \itemindent
   \itemsep 0pt
   \parsep 0pt}
   \def\newblock{}
   \sloppy
   \sfcode`\.=1000\relax}


\setlength{\baselineskip}{12.7pt}

\title{
Public Access Defibrillator Deployment for Cardiac Arrests: A Learn-Then-Optimize Approach with SHAP-based Interpretable Analytics
}

\author{\begin{center}Kexin Cao\textsuperscript{1}\footnotemark[2]\ , 
    Chih-Yuan Yang\textsuperscript{1}\footnotemark[2]\ , 
    Keng-Hou Leong\textsuperscript{1,2}, 
    Xinglu Liu\textsuperscript{1,3}, 
    Wai Kin (Victor) Chan\textsuperscript{1}\footnotemark[1]\\
[11pt]
\textsuperscript{1}Institute of Data and Information, Tsinghua Shenzhen International Graduate School, Tsinghua University,  Shenzhen, China\\
\textsuperscript{2}Institute of Remote Sensing and Geographical Information Systems, School of Earth and Space Sciences, Peking University, Beijing, China\\
\textsuperscript{3}Department of Civil Engineering, The University of Hong Kong, Hong Kong, China \end{center}}


\footnotetext[2]{These authors contribute equally to this work.}
\footnotetext[1]{Corresponding author.}

\maketitle

\vspace{-12pt}

\section*{ABSTRACT}
\noindent
Out-of-hospital cardiac arrest (OHCA) survival rates remain extremely low due to challenges in the timely accessibility of medical devices. Therefore, effective deployment of automated external defibrillators (AED) can significantly increase survival rates. Precise and interpretable predictions of OHCA occurrences provide a solid foundation for efficient and robust AED deployment optimization. This study develops a novel learn-then-optimize approach, integrating three key components: a machine learning prediction model, SHAP-based interpretable analytics, and a SHAP-guided integer programming (SIP) model. The machine learning model is trained utilizing only geographic data as inputs to overcome data availability obstacles, and its strong predictive performance validates the feasibility of interpretation. Furthermore, the SHAP model elaborates on the contribution of each geographic feature to the OHCA occurrences. Finally, an integer programming model is formulated for optimizing AED deployment, incorporating SHAP-weighted OHCA densities. Various numerical experiments are conducted across different settings. Based on comparative and sensitive analysis, the optimization effect of our approach is verified and valuable insights are derived to provide substantial support for theoretical extension and practical implementation.

\section{INTRODUCTION}
Out-of-Hospital Cardiac Arrest (OHCA) is a significant public health concern and remains one of the leading causes of death worldwide \shortcite{Yan2020}. 
\replaced[id=Cao]{There were 35,801 OHCA observations treated by the Emergency Medical System (EMS) over 62.11  million person-years, resulting in an incidence rate of 54.99 per 100,000 person-years \shortcite{rea2004incidence}.}{In China, an estimated 700,000 sudden cardiac deaths occur annually, with 87\% taking place outside of medical facilities \shortcite{NCCDChina2024}.}
Unfortunately, the survival rates for OHCA patients \deleted[id=Cao]{in China }are alarmingly low at just 1.2\% \shortcite{Zheng2023}.
Timely medical interventions, such as \deleted[id=Cao]{the use of }Automated External Defibrillators (AEDs), have been shown to significantly improve survival rates, potentially raising them to 16\% or higher \shortcite{Nichol2008}.
\replaced[id=Cao]{However, achieving such objectives in practice presents significant challenges}{achieving these outcomes in practice is challenging}, as the critical “golden four minutes” \added[id=Cao]{time }window for effective treatment \replaced[id=Cao]{is frequently unsatisfied, particularly in regions with insufficient AED coverage}{is often exceeded, especially in areas where AED coverage is insufficient} \shortcite{Berdowski2010}.
These challenges highlight the urgent need to improve AED accessibility\replaced[id=Cao]{, thereby reducing AED delivery}{and reduce response} times to increase OHCA survival rates.

The effectiveness of medical intervention is closely tied to the proximity of AEDs to the OHCA incident spots. \replaced[id=Cao]{Studies have demonstrated that deploying AEDs in and around high-risk areas significantly increases the likelihood of a successful rescue \shortcite{blom2014improved,lin2023optimal}.}{ show that placing AEDs in close proximity to high-risk areas can greatly increase the likelihood of successful resuscitation }
Therefore, \replaced[id=Cao]{by analyzing the spatial distribution of OHCA and identifying high-risk areas, the coverage efficiency of AED deployment can be significantly enhanced.}{optimizing AED deployment by predicting the spatial distribution of high-risk OHCA areas becomes crucial.}
\replaced[id=Cao]{The current primary methods for identifying high-risk OHCA areas rely heavily on demographic data}{Traditional approaches to identifying high-risk OHCA locations rely heavily on population data, such as population density and demographic structure,} as inputs for machine learning models or on statistical analysis of historical OHCA occurrence data \shortcite{hessulf2023predicting,sasson2010predictors,nakashima2023machine}.
\replaced[id=Cao]{While previous studies have demonstrated the effectiveness of these predictive models, their high demand for demographic or historical data imposes limitations on their practical extension to regions where such data are scarce.}{While effective in certain researches, these methods have notable limitations.} 
\replaced[id=Cao]{In reality, the collection of demographic data is inherently time-consuming and laborious due to its complicated processes. Besides these difficulties, the recording of historical OHCA data faces additional challenges, such as high needs for information technology and hardware facilities in cities, necessitated by the very short time window and the high timeliness and accuracy required for recording.}{Collecting population and demographic data is both labor-intensive and time-consuming.
Additionally, the lack of comprehensive historical OHCA datasets in many cities presents significant barriers to accurate risk prediction.}
\added[id=Cao]{To achieve accurate OHCA risk prediction in cities lacking these data, this study attempts to employ more accessible geographic feature data, i.e., Point-of-Interest (POI) and building distribution data, to potentially represent demographic characteristics and facilitate the development of a practical and scalable OHCA risk prediction framework. Such a framework enables our trained model to achieve domain generalization in any alternative region, especially without demographic features and historical OHCA records.}

\added[id=Cao]{To clarify the relationship between geographical features and OHCA occurrence, we apply machine learning (ML) methods due to their generalizability and excellent predictive performance in urban planning issues \shortcite{nosratabadi2019state,chen2024optimized}. Furthermore, for such a life-critical problem, we aim to interpret the ML model, determining the contribution of each input feature to the prediction output, which is challenging when using solely 'black box' machine learning models.}
The SHAP (Shapley Additive Explanations) framework \shortcite{lundberg2017unified}, grounded in the classical Shapley value from game theory, offers a \replaced[id=Liang]{viable}{promising} solution to these interpretability challenges.
It helps bridge the gap between the ML prediction and the decision-making framework, even in complex environmental contexts \shortcite{ekmekcioglu2022explainable,iban2022machine,CHEN2024105191}.
By quantifying the contribution of \replaced[id=Cao]{each geographic feature to OHCA risk density, SHAP values provide a transparent and evidence-based basis for making interpretable and accurate AED deployment decisions. This improves the reliability and practical feasibility of optimizing AED deployment strategies \shortcite{Ribeiro2016}.}{each feature to the model's output, SHAP values offer valuable insights into the relationship between geographic features and OHCA risk density. This interpretability provides a transparent and evidence-based foundation for making interpretable and accurate AED deployment decisions, improving both the reliability and practical feasibility of optimized AED placement strategies \shortcite{Ribeiro2016}.}

\begin{figure}[htb]
  \centering
  \includegraphics[width=0.8\textwidth]{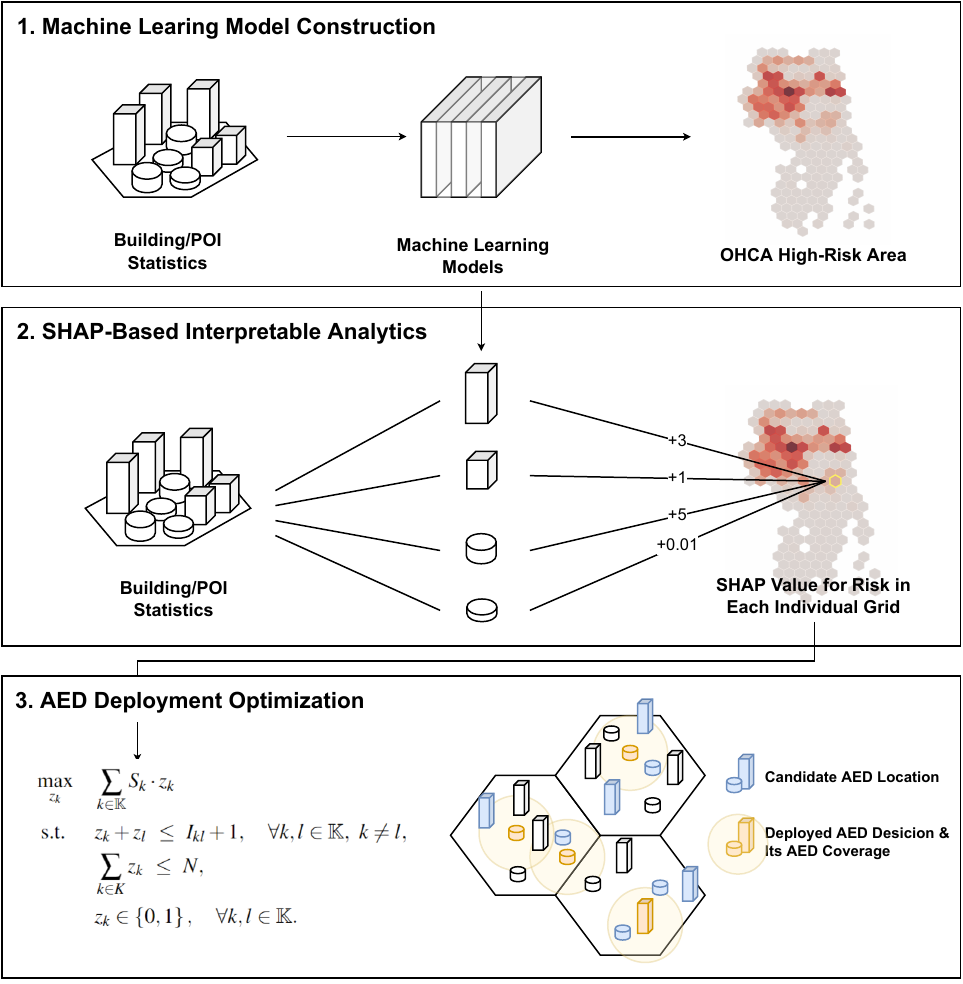}
  \caption{Learn-then-optimize framework with SHAP-based interpretable analytics}
  \label{fig:system}
\end{figure}

\replaced[id=Cao]{In summary, this study aims to develop a learn-then-optimize framework that combines machine learning, SHAP-based interpretable analytics, and integer programming (see Fig.\ref{fig:system}). The framework predicts the risk density of OHCAs using POI and building distribution data and calculates SHAP values as input parameters for optimizing AED deployment decisions. The primary aim of this study is to enhance the efficiency of AED coverage and to increase patient survival rates through the following three key components of contribution: prediction, interpretation, and optimization.}{The objective of this study is to develop an integrated machine learning and optimization framework for AED deployment using only the Points of Interest (POI) and building distribution data. This framework aims to enhance accessibility to AED and maximize OHCA survival rates through three key contributions:}
\begin{enumerate}
    \item \replaced[id=Cao]{We confirm the feasibility of predicting OHCA high-risk areas using a Neural Network (NN) model trained only on POI and building distribution data as input. The predictive performance of our NN model on the test set is also validated with an $R^2$ value exceeding 0.75. Moreover, our findings reveal a strong correlation between geographic features and OHCA events, indicating that it is compelling to further construct an interpretable SHAP model with geographic features based on the NN model.}{We confirm the feasibility of using an NN model trained solely with POI and building distribution data as input to predict high-risk areas for OHCA. The prediction performance of our NN model on the test set has also been verified with an $R^2$ value over 0.75. Our findings demonstrate a strong correlation between geographic features and OHCA incidents, thus prove that it is meaningful to further construct an interpretable SHAP model based on the NN model trained with geographic feature data.}

    \item \replaced[id=Cao]{An interpretable SHAP model is developed to quantify the contribution of various geographic features to the incidence of OHCA in each area. These SHAP values for different building types are utilized as input parameters for the AED deployment model to enhance its practical feasibility and optimization effectiveness. Additionally, further analysis of SHAP values yields practical insights, aiding decision-makers in identifying critical building types and thus high-risk areas for urban emergency deployment.}{We construct an interpretable model for our NN using SHAP values, which quantify the relationship of different geographic features to the OHCA occurrence in each area. The SHAP values offer valuable insights that can be applied to an operational model for optimal AED deployment.}
    
    \item \replaced[id=Cao]{A SHAP-guided integer programming (SIP) model is established to optimize the AED deployment, and experiments are conducted on different deployment scales. Based on the experimental results, we compare the deployment effects, including historical OHCA coverage efficiency and average survival rates, between the SIP model and a random deployment baseline to validate the optimization effect of our SIP model. Furthermore, a sensitivity analysis of the minimum spacing for deployment is employed, providing a theoretical reference for setting the minimum spacing of AEDs during deployment. Some practical insights for policy-makers to implement efficient and cost-effective AED deployment strategies are also obtained.}{We build a SHAP-value-guided Integer Programming (SIP) model to optimize the deployment of AEDs. Based on the simulation results, we compare the deployment effects including historical OHCA coverage rates and average survival rates between our SIP model and the random-deployment baseline to validate the optimization effect of our model.}
\end{enumerate}

Overall, our proposed learn-then-optimize framework strives to overcome the existing challenges in AED deployment, provide accurate and practical deployment decisions, and enhance urban healthcare emergency service.
The rest of this paper is organized as follows. Section 2 presents the methodology of our study, including an NN prediction model, a SHAP interpretation model, and the SIP AED deployment model. Section 3 examines the experimental results along with the insights derived from each component. Finally, Section 4 summarizes the study and suggests potential directions for future work.

\section{methodology}

This subsection elaborates on the three components of our proposed learn-then-optimize framework. First, we identify OHCA high-risk areas using an NN model trained on geographic data as input. Second, we leverage SHAP values to quantify the impact of geographic features and building types on OHCA risk predictions, ensuring transparency and interpretability. Finally, we incorporate these SHAP-derived insights into an integer programming (IP) optimization model for AED placement. The performance of the proposed model is then compared with random AED distributions under various scenarios.

\subsection{Data Preprocessing (Geographic Information System)}

To predict OHCA high-risk areas \replaced[id=Cao]{rationally and accurately}{in a rational and accurate manner}, we apply the H3 grid system developed by Uber Technologies \shortcite{brodsky2018h3} 
\added[id=Cao]{ as a statistical unit for OHCA event occurrence}. This grid system offers \replaced[id=Cao]{the advantages of standardization, divisibility, and flexibility.}{several meaningful advantages:} 
\begin{enumerate}
    \item The hexagonal grid design ensures uniform grid areas, enabling an accurate quantification and comparison of OHCA occurrences across various grids. The hexagon grids \replaced[id=Cao]{can}{could} be approximated as circles as the radii of their circumcircles and inscribed circles are close.
    \item The system is divided into 16 levels of resolution, with each parent grid subdivisible into seven smaller child grids. The levels range from the largest grid size at level 0 to the smallest at level 15. This hierarchical structure \replaced[id=Cao]{enables the predictions of}{allows us to predict} OHCA risk \replaced[id=Cao]{at both macro and}{on either macro or} micro scales, offering flexibility in analysis. 
    \item Compared to traditional administrative \replaced[id=Cao]{divisions}{boundaries} (\textit{e.g.}, township\deleted[id=Cao]{ divisions}), the H3 grid system more accurately reflects \replaced[id=Cao]{the cross-boundary AED coverage,}{AED coverage across boundaries, including those} spanning multiple administrative areas. 
\end{enumerate}

The primary objective of our project is to optimize AED deployment and reduce delays in medical intervention. To ensure the grid size appropriately aligns with the coverage range of AEDs, we select grid size level 7 with an average edge length of 1.41 km, which allows \replaced[id=Cao]{an emergency responder}{a volunteer} to run from the center to the edge in 4 minutes.
Research has shown that this choice strikes a balance between granularity and practicality, making our model effective for addressing real-world emergency scenarios \shortcite{Zheng2023}.

\subsection{Machine Learning Model Construction}

We obtain the geographic features from OpenStreetMap \shortcite{haklay2008openstreetmap} as the input of the model, which includes types and locations (latitude and longitude) of most POIs and buildings in Virginia Beach and other cities.
The dataset consists of 76 types of POIs on the map, including restaurants, community centers, post offices, etc, and 39 types of buildings, such as apartments, schools, and hotels.
The list of the types of POIs and buildings is provided in Appendix \ref{appa}.
For the model output, we apply the dataset of OHCA occurrences from January 2017 to June 2019 in Virginia Beach \shortcite{custodio2022spatiotemporal}.

The set of grids is represented as $\mathbb{I}=\{0,1,\cdots,G\}$, where $\mathbb{I}_i$ denotes each grid cell $i$. 
Define $X_i$ as the counts of various geographic features within grid $i$, and $y_i$ as the historical OHCA occurrences in the corresponding grid.
To further elaborate, $y_i$ originally denotes the OHCA density of grid $i$, which is simplified to the total number of OHCAs in the following prediction model, as all grids have the same area.
With a hypothesis that the OHCA occurrence of a grid is solely contributed by the geographical features within the grid, we can formulate
the predictive NN model $f$, 
\begin{equation}
y_i=f(X_i),
\end{equation}

We employ a Multilayer Perceptron (MLP) to perform a regression on the relationships between $X_i$ and $Y_i$.
The trained MLP is then analyzed and explained by a SHAP model.


\subsection{SHAP-Based Interpretable Analytics\label{sec2_3}}
SHAP \shortcite{lundberg2017unified} offers a reliable framework for interpreting model predictions by assigning each feature a SHAP value, which quantitatively represents its contribution to the model's output. 
\deleted[id=Liang]{The SHAP method satisfies an additive property, which ensures that the sum of all SHAP values for a given predicted value, along with a baseline prediction, exactly reconstructs the model's predicted value.}
The SHAP reconstruction $g$ can be expressed as:

\begin{equation}
f(X_i) \approx g(X_i) = \phi_{0} + \sum_{j=1}^M \phi_{ij},
\label{eq:shap_additive}
\end{equation}
where
\begin{itemize}
    \item \(f(X_i)\): The model's output (prediction) for a given input \(X_i\);
    \item \(\phi_0\): The baseline prediction, typically the average of all predictions across the dataset;
    \item \(\phi_{ij}\): The SHAP value for feature \(j\) in grid $i$, representing the contribution of that feature to the prediction \(f(X_i)\);
    \item \(M\): The total number of features in the model.
\end{itemize}
A significant advantage of SHAP lies in its additive \replaced[id=Liang]{property, which ensures that the sum of all SHAP values for a given predicted value, along with a baseline prediction, exactly reconstructs the model's predicted value.}{nature, which guarantees that the sum of all SHAP values equals the difference between the model's prediction and the average (baseline) prediction.}
This property holds true even for the most complex machine learning models, allowing SHAP to consistently explain predictions in a transparent and interpretable way.
By breaking down the model's output into the contributions of individual features, SHAP provides a clear and reliable understanding of how each feature \added[id=Liang]{additively} influences the final prediction. 

Due to the massive computation, we apply an approximation method \shortcite{sundararajan2020many} to compute $\phi_{ij}$, resulting in the approximation in Equation \ref{eq:shap_additive}.
The calculation of SHAP values is based on cooperative game theory.
The formula is:
\begin{equation}
\phi_{ij} = \sum_{S \subseteq F \setminus \{j\}} \frac{|S|!(M - |S| - 1)!}{M!} \left[ f(X_{i, S \cup \{j\}}) - f(X_{iS}) \right],
\label{eq:shap_value}
\end{equation}
where
\begin{itemize}
    \item \(F\) is the set of all features in the model;
    \item \(S \subseteq F \setminus \{j\}\) is a subset of the set \(F\) (all features) excluding feature \(j\), representing all possible subsets without \(j\);
    \item \(f(X_{i, S \cup \{j\}})\) is the model output with a modified input $X_{i, S \cup \{j\}}$, in which the feature subset \(S\) and feature \(j\) in $X_i$ are retained, while the other features are substituted with white noises;
    \item \(f(X_{iS})\) is the model output with a modified input $X_{iS}$, in which \(S\) are retained while the other features (including $j$) are substituted with white noises;
    \item \(|S|!(M - |S| - 1)!/M!\) is the combinatorial factor that accounts for the possible feature combinations when calculating the SHAP value. This weight ensures the fair distribution of contributions among features.
\end{itemize}

The calculation of SHAP values considers all possible combinations of features, both with and without each feature. This approach comprehensively evaluates the interactions and dependencies between features, ensuring their collective influence on the prediction is fully captured. Additionally, SHAP guarantees a fair distribution of feature contributions, accurately reflecting their impact on the model’s output. Therefore, if a feature has a negative impact on the prediction, decreasing the predicted value, its corresponding SHAP value will be negative, and vice versa \shortcite{lundberg2017unified}.

In our study, we calculate the SHAP values for various building features using the SHAP library in Python. The SHAP values allows us to quantify the positive or negative contribution of each feature to the prediction \(y_{ij}\). This enables a numerical analysis of the relationship between each building feature and the OHCA high-risk area. Additionally, the SHAP values can serve as the basis for AED site selection. By distributing the SHAP values across buildings, we assign these values as scores to each individual building. These scores reflect how much an individual building contributes to the predicted number of OHCA cases in the grid. In other words, buildings with higher scores have a greater influence on the predicted OHCA occurrences in the area.
 

\subsection{AED Deployment Optimaization\label{sec2_4}}

This subsection first describes the AED deployment optimization problem based on predicted SHAP values and then formulates a SHAP-guided integer programming (SIP) model with its objective function and various constraints elucidated.

\begin{figure}[htb]
  \hspace{1.7em}
  \includegraphics[width=0.95\textwidth]{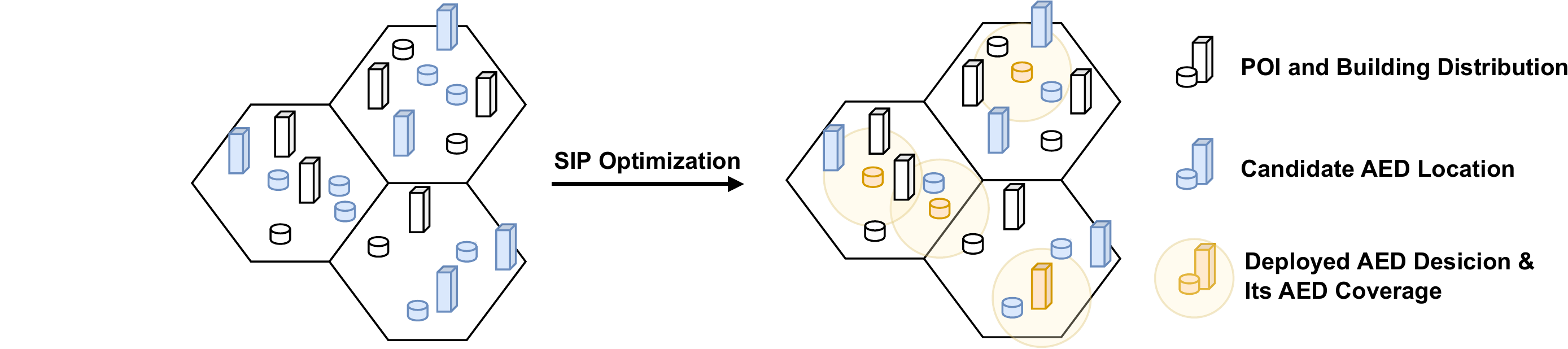}
  \caption{A toy example for the AED deployment optimization problem decisions}
  \label{fig:prob}
\end{figure}

Fig.\ref{fig:prob} illustrates a toy example for the AED deployment optimization problem decisions. A summary of all sets, parameters, and variables related to the SIP model is listed in Table \ref{table_set}. To begin with, two hypotheses are proposed for computational convenience while maintaining the plausibility of the problem:
\begin{enumerate}
\item \replaced[id=Liang]{The SHAP values $\phi _{ij}$ in Section \ref{sec2_3} are computed based on $X_{ij}$, which is the quantity of each type of geographical feature within a grid, for instance, three hospitals contributing to six OHCAs. From the linear additivity of SHAP calculations, we further assume that the SHAP value of any geographical feature type $i$ in a certain grid $j$ can be equally shared across every single building $p$ of that type (see Equation \ref{phip}), for example, one hospital contributing to two OHCAs;
}{The SHAP value we derived from the NN ideally reflects the practical circumstance and could be shared across all the buildings or POIs of the same type within the grid.}
\begin{align}
    \phi _p\ &=\ 
    \begin{cases}
	\frac{\phi _{ij}}{X_{ij}},\ &\text{if\ }X_{ij}>0,\\
	0,\ \ &\text{if\ }X_{ij}=0,\\
    \end{cases} 
    \label{phip}
\end{align}
\item In accordance with the exclusive hypothesis for each geographical feature relative to its respective grid, only buildings within the AED coverage area are taken into account when calculating the occurrence of OHCA.
\end{enumerate}

Among the urban planning areas, there are numerous POIs and buildings, from which candidate AED locations $\mathbb{K}$ are randomly selected. In line with Hypothesis 2, for each candidate location $k \in \mathbb{K}$, we use SHAP values $\phi_{p}$ of the POIs and buildings within a certain range $A_k$ to calculate its SHAP-weighted OHCA density $S_k$ (See Equation \ref{shap_weighted}), which represents the weighted coverage of predicted OHCAs.
\begin{align}
    S_k\ &=\ \frac{\sum_{\mathbb{I}_i\subseteq \mathbb{I},\ \mathbb{I}_i\cap \mathbb{A}_k\ne \emptyset}{\left( |\mathbb{I}_i\cap \mathbb{A}_k|\cdot \frac{\sum_{p\subseteq \mathbb{I}_i\cap \mathbb{A}_k}{\phi _p}}{|\mathbb{I}_i|} \right)}}{|\mathbb{A}_k|} \nonumber\\ 
    &= \ \sum_{\mathbb{I}_i\subseteq \mathbb{I},\ \mathbb{I}_i\cap \mathbb{A}_k\ne \emptyset}{\left( \frac{|\mathbb{I}_i\cap \mathbb{A}_k|}{|\mathbb{A}_k|}\cdot \frac{\sum_{p\subseteq \mathbb{I}_i\cap \mathbb{A}_k}{\phi _p}}{|\mathbb{I}_i|} \right)}\nonumber\\ 
    &=\ \sum_{\mathbb{I}_i\subseteq \mathbb{I},\ \mathbb{I}_i\cap \mathbb{A}_k\ne \emptyset}\omega_{ik} \cdot \rho_{ik},
    \label{shap_weighted}
\end{align}
where $|\mathbb{I}_i\cap \mathbb{A}_k|$ denotes the area of the overlapping region between grid $\mathbb{I}_i$ and $A_k$. $\rho_{ik}$ represents the OHCA density for each overlapping region and $\omega_{ik}$ can be interpreted as an equivalent area weighting factor.
\begin{table}[htbp]
\label{table1}
\renewcommand\arraystretch{1.3}
  \centering
   \caption{Sets, parameters, and variables}
    \begin{tabular}{p{1cm}<{\centering} p{11cm}}
    \toprule
    \multicolumn{2}{l}{\textbf{Sets}} \\
     $\mathbb{P}$ & Set of POI and building locations, $\mathbb{P}=\{0,1,\cdots,P\}$\\
     $\mathbb{K}$ & Set of AED candidate locations, $\mathbb{K}=\{0,1,\cdots,K\}$\\
     $\mathbb{H}$ & Set of historical OHCA incidents, $\mathbb{H}=\{0,1,\cdots,H\}$\\
    \toprule
    \multicolumn{2}{l}{\textbf{Parameters}} \\
     $P$ & Number of POI and building locations \\
     $K$ & Number of AED candidate locations \\
     $H$ & Number of historical OHCA incidents \\
     $N$ & Total number of AEDs to deploy \\
     $D_{min}$ & Minimum spacing restriction for two AED locations \\ 
     $X_{ij}$ & Number of each geographic feature $j$ in each grid $i$ \\
     $\phi_{p}$ & SHAP value of each POI and building location $p \in \mathbb{P}$ \\
     $S_k$ & SHAP-weighted OHCA density of each AED candidate location $k \in \mathbb{K}$ \\
     $D_{kl}$ & Distance between two AED locations $k$ and $l$, $k,l \in \mathbb{K}$, $k \ne l$\\
     $I_{kl}$  & Indicator parameter, equals to 1 if $D_{kl}$ is less than $D_{min}$, otherwise 0\\   
     \toprule
    \multicolumn{2}{l}{\textbf{Decision Variables}} \\
     $z_k$ & Equals to 1 if AED candidate location $k \in \mathbb{K}$ is chosen, otherwise 0\\
    \bottomrule
    \end{tabular}%
  \label{table_set}%
\end{table}%




Based on the computed SHAP-weighted OHCA densities, an AED deployment optimization model is formulated with the objective of maximizing the weighted coverage of predicted OHCAs. Constraints (\ref{distancelimit}) restrict the minimum spacing between each two AED deployed locations and Constraints (\ref{locsum}) determine the total number of AEDs to deploy. The binary integer constraints for the decision variables are given in Constraints (\ref{var01}). 
\begin{align}
    \underset{z_k}{\max}\quad &\sum_{k\in \mathbb{K}}{S_k\cdot z_k}
    \label{obj}\\
    \mbox{s.t.   } \quad  &z_k+z_l\ \le \ I_{kl}+1, \quad \forall k,l\in \mathbb{K},\ k\ne l,
    \label{distancelimit}\\
    &\sum_{k\in K}{z_k\ \le \ N,} 
    \label{locsum}\\
    &z_k\in \left\{ 0,1 \right\}, \quad \forall k,l\in \mathbb{K}.
    \label{var01}
\end{align}




To further examine the model optimization effect in numerical experiments, we compute the historical OHCA coverage and the average survival rate as performance measures. Equation (\ref{coverage}) represents the historical OHCA coverage, which is the total number of OHCAs covered by at least one deployed AED within a given coverage radius $C_{R}$. Survival rates are computed based on the simulated AED delivery time from deployed AEDs to OHCAs according to Equation (\ref{eq:survival}).

\begin{equation}
\text{Historical\ OHCA\ coverage} = \sum_{h \in \mathbb{H}} \mathbf{1} \left ( \min_{k \in \mathbb{K}} \left (Dist\{k,h\}|z_k = 1 \right) \le C_{R} \right )
\label{coverage}
\end{equation}



\begin{equation}
\text{Survival\ rate} = 
\begin{cases}
  (1 + e^{-0.26 + 0.106 \cdot t_{\text{AED}} + 0.139 \cdot t_{\text{CPR}}})^{-1}, & t_{\min} < 4, \\
    0, & t_{\min} \geq 4.
  \label{eq:survival}
\end{cases}
\end{equation}

\section{Experimental Results and Analysis}
\subsection{Data Source}
To train and test an NN model with only POI and building distribution data to predict the high-risk area of OHCA, a regional dataset that includes both the history of OHCA occurrences and geographic features is necessary. 
We obtain both datasets from Virginia Beach, a coastal city located in Virginia, U.S..
Virginia Beach is the largest city in Virginia, covering an area of approximately 497 square miles (1,287 square kilometers) and comprising a grid of 177 H3 class 7 hexagons.
With a population of around 460,000, it is also the most populous city in the state.

\subsubsection{Out-of-Hospital Cardiac Arrest Data}

Our research utilizes a geographical dataset from Virginia Beach \shortcite{custodio2022spatiotemporal}, which records the locations of all OHCA cases reported between January 1, 2017, and June 30, 2019. To prepare these data for analysis, we assign a geographic label to each OHCA case based on its latitude and longitude, using the H3 Level 7 hexagonal grid system. OHCA cases within the same H3 Level 7 hexagonal grid share the same label. Subsequently, we sum up the OHCA cases within each grid cell, with the total count in each cell serving as the target variable for our neural network model. This transformation ensures standardization, divisibility, and flexibility, as highlighted in Section 2.1, thereby enhancing our model's capability to effectively address real-world emergency scenarios.

\subsubsection{POI and Building Distribution Data}
The POI and building distribution data used in our research are extracted from OpenStreetMap (OSM), a collaborative, open-source mapping platform that provides detailed and up-to-date geographic data contributed by a global community of users. OSM offers comprehensive information on buildings, landmarks, and other spatial features, making it a valuable resource for spatial analysis.

We collect building information for the entire Virginia Beach area, including 115 buildings or POI features such as apartments, schools, fountains, fuel stations, and more. The full feature description can be checked on \href{https://wiki.openstreetmap.org/wiki/Map\_features\#Public\_transport}{https://wiki.openstreetmap.org/wiki/Map\_features\#Public\_transport}. To ensure consistency in the experiment and to leverage the advantages of the H3 level 7 grid system, we apply the same geocoding system to process the geographic data. The total count of each building feature within a grid cell serves as input for our NN model.

\subsection{OHCA High-Risk Area Prediction}
To validate the feasibility of utilizing an NN model based solely on POI and building distribution data to predict high-risk areas for OHCA, this study aims to establish a significant correlation between geographic features and OHCA risk areas.

Given the development status of Virginia Beach, we employ the eastern portion of Virginia Beach as the training set, comprising 83 grids, and the western part as the testing set, including 94 grids, for our NN model training. This division ensures that adjacent grids are included in both the training and testing sets, accounting for the fact that AED coverage often spans multiple grids. Consequently, this method facilitates more practical and accurate deployment of AEDs in real-world scenarios.

\begin{figure}[ht!]
  \centering
  \includegraphics[width=0.7\textwidth]{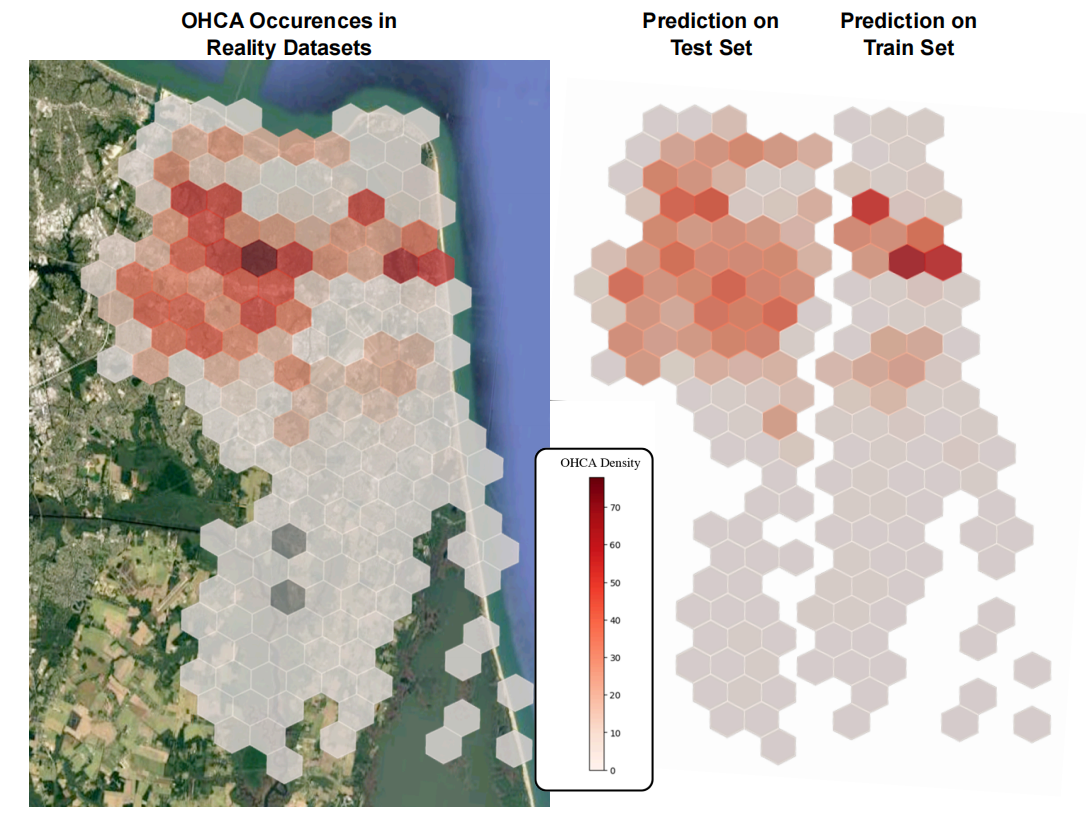}
  \caption{OHCA occurrences and and prediction results in each H3 Level 7 grid in Virginia Beach.}
  \label{fig:test_vs_train}
\end{figure}

The result and visualization of our three-layer fully connected NN model are presented in Fig. \ref{fig:test_vs_train} and Table \ref{tab:machine_learning}. 
The training set results demonstrate a strong fit ($R^2$ = 0.975, MAE = 1.28), indicating that the model has learned the training data well.
Additionally, the test set results ($R^2$ = 0.752, MAE = 5.56) suggest a reasonable correlation between the predictions and actual values.
This indicates that the NN model effectively captures the relationship between geographic data and OHCA occurrences.
MAPE (Mean Absolute Percentage Error) represents the error as a percentage by averaging the absolute percentage differences between predicted and actual values.
However, considering that our dataset may contain zero values, we opt to use MAE instead of MAPE to assess the model's performance.

\begin{table}[htbp]
\renewcommand\arraystretch{1.3}
  \centering
   \caption{Model performance on train and test sets}
    \begin{tabular}{p{3cm} p{1.8Cm}<{\centering} p{1.8cm}<{\centering}}
    \toprule
    Dataset         & $R^2$         & MAE \\
    \toprule
    Train set   & 0.975         & 1.28 \\ 
    Test set    & 0.752         & 5.56 \\
    \toprule
    \end{tabular}%
  \label{tab:machine_learning}%
\end{table}%

The results confirm the feasibility of employing an NN model that exclusively uses POI and building distribution data to predict high-risk areas for OHCA incidents. This study substantiates the presence of a significant correlation between geographic features and OHCA occurrences. Furthermore, this research proves that constructing an interpretable model based solely on geographic feature data through an NN framework is highly valuable.

\subsection{Model Interpretation}

\begin{figure}[ht!]
  \centering
  \includegraphics[width=0.75\textwidth]{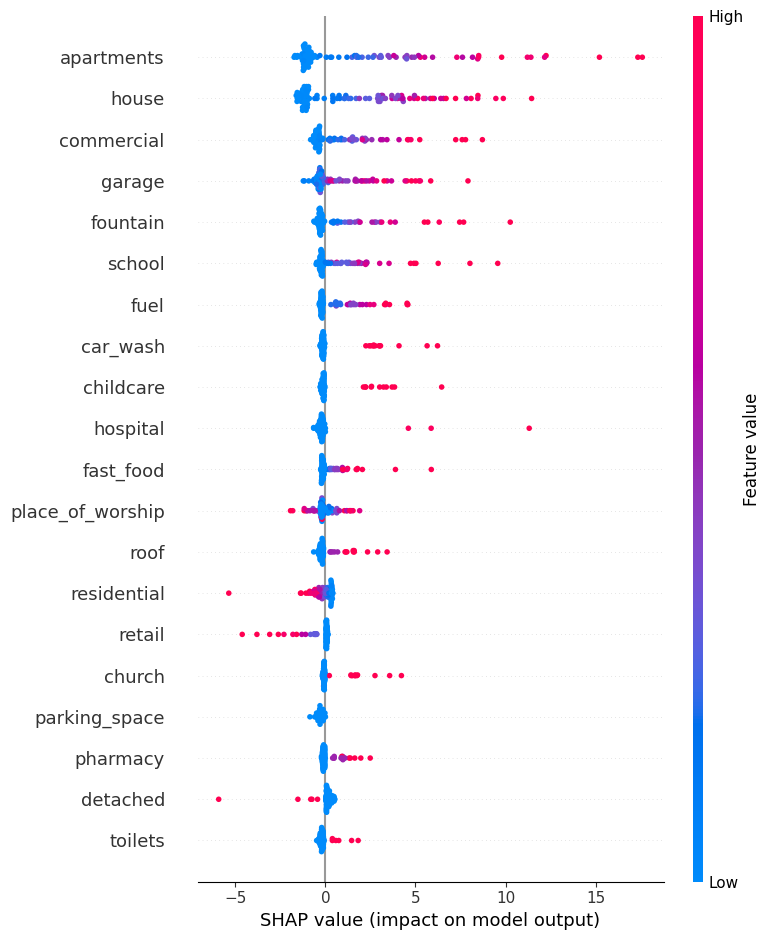}
  \caption{Top 20 geographic features with the highest influence ranked by absolute SHAP values.}
  \label{fig:shap}
\end{figure}

Fig.\ref{fig:shap} illustrates the SHAP values of the most influential features on the prediction.
Detailed SHAP values for all features are provided in Appendix \ref{appb}. As the mean absolute SHAP value quantitatively represents the influence of a feature on the prediction \cite{ekanayake2022novel}, the features in this figure are ranked according to their absolute SHAP values.
This ranking highlights the most critical determinants that shape the predictive outcomes.

The colored points, ranging from blue to red, indicate the counts of POI or building features within a grid, with blue representing lower counts and red representing higher counts.
The x-axis displays the SHAP value, which quantifies each feature's contribution to the model's prediction for a specific grid.
\replaced[id=Liang]{Positive SHAP values indicate that the feature tends to increase the predicted number of OHCA occurrences, while negative values suggest a decrease.}{Positive SHAP values indicate that the feature drives the prediction toward a higher number of OHCA occurrences, while negative values drive the prediction toward a lower number of OHCA occurrences.}
\deleted[id=Liang]{In this experiment, we apply SHAP values for model interpretability analysis to reveal the contribution of different features to the prediction results. The visualization of SHAP values is used to demonstrate the impact of each feature on the prediction outcome.} 
For example, points with high counts of \textit{apartment} are associated with high SHAP values, \deleted[id=Liang]{, associated with a larger number of OHCA occurrences.}
indicating a positive correlation between the features \textit{apartment} and OHCA occurrence: areas with higher \textit{apartment} density are more likely to have a greater number of OHCA occurrences.
Conversely, a higher count of the features \textit{retail} is associated with \replaced[id=Liang]{negative SHAP values}{a lower number of OHCA cases}, indicating a negative correlation between \textit{retail} density and OHCA incidents. 
\deleted[id=Liang]{, meaning areas with higher retail density are more likely to experience fewer OHCA cases}

The results validate the effectiveness of geographic data in capturing critical insights into population demographics and mobility dynamics, which are valuable for analyzing OHCA incidents.
As illustrated in Fig.\ref{fig:shap}, we have identified the leading critical determinants that shape the predictive outcomes.
The findings indicate a positive correlation between population density and OHCA incidents, particularly in regions dominated by apartments and residential housing.
In contrast, open and less densely populated areas, such as graveyards and retail spaces, are typically associated with fewer OHCA incidents.
In the United States, retail areas are often located in commercial districts, which generally have lower residential density compared to residential neighborhoods.
This observation leads us to conclude that areas with higher residential density, such as residential neighborhoods, are high-risk areas for OHCA, while areas with lower residential density, such as those with parking lots, are low-risk areas for OHCA.
This conclusion aligns with findings suggesting that regions with higher population density tend to experience a greater number of OHCA cases \shortcite{sasson2010predictors}.
According to the above-mentioned findings, the implied connection between geographic and demographic features in predicting OHCA risk densities can inform AED deployment decisions.


\subsection{AED Deployment}

To validate the robustness and efficacy of the proposed SIP model, we have conducted a series of numerical experiments across various settings, compared with the random baseline.
The experiments cover various deployment scales $N$ (ranging from 5 to 100) and minimum spacing $D_{min}$ (ranging from 0 to 1.6 km), designed to conduct a sensitive analysis and verify the model's superiority.
For each case, we select 10 sets of location candidates, each with 5000 candidates randomly sampled from a total of 99724 POIs and buildings in Virginia Beach.
For the random baseline, the solutions are generated by randomly selecting $N$ candidates, without imposing any restrictions on the minimum spacings.

\begin{figure}[htb]
  \centering
  \includegraphics[width=0.8\textwidth]{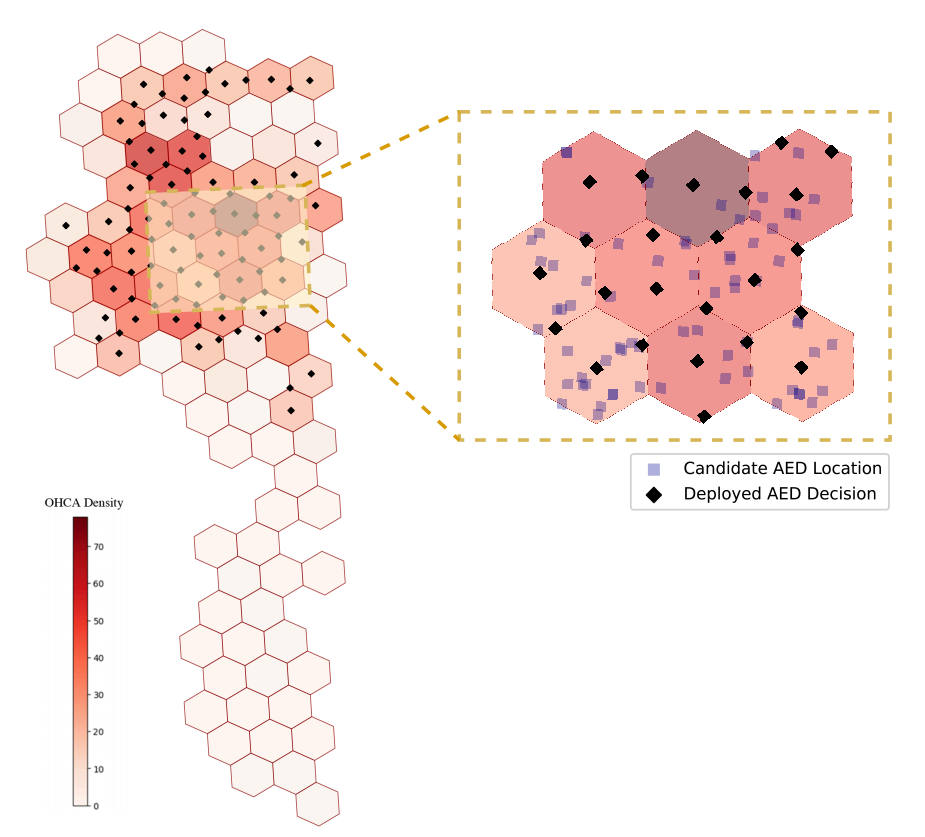}
  \caption{ AED deployment decisions based on the SIP model with $N=100$ and $D_{min}$ = 1.2 km.}
  \label{fig:visualization}
\end{figure}

\subsubsection{Visualization of Deployment Decisions}

Fig.\ref{fig:visualization} visualizes the solution of the SIP Optimization when $N=100$ AEDs are deployed with a minimum spacing of $D_{min}$ = 1.2 km.
The AEDs deployed from our Learn-then-Optimize framework cover most high-SHAP-value geographical features and align well with real-world high-risk areas, following a nearly hyper-uniform pattern.
This results supports the hypothesis presented in Section \ref{sec2_4}, in which the OHCA occurrence in a grid can be shared among the geographical features within that grid.

We observe that with $N=100$ AEDs deployed, at least one AED deployed is present in all high-risk grids (where $y_i > 16$), covering an average of 1388.2 out of 1501 cases (92.4\%) across the western half of Virginia Beach.
The effective performance achieved by the SIP solutions, guided by the SHAP values, validates the analytics discussed in Section \ref{sec2_3}.
However, there remains a gap of approximately 110 cases that could potentially be saved with additional AEDs.
As provided in the following results, we consider $N$ = 100 as a critical point for saturation, beyond which the marginal benefit of deploying more AEDs in low-risk grids (southeast of Virginia Beach) is little.

\subsubsection{Optimization Effect: SIP v.s. Random Baseline}

\begin{figure}[ht!]
  \centering
  \includegraphics[width=\textwidth]{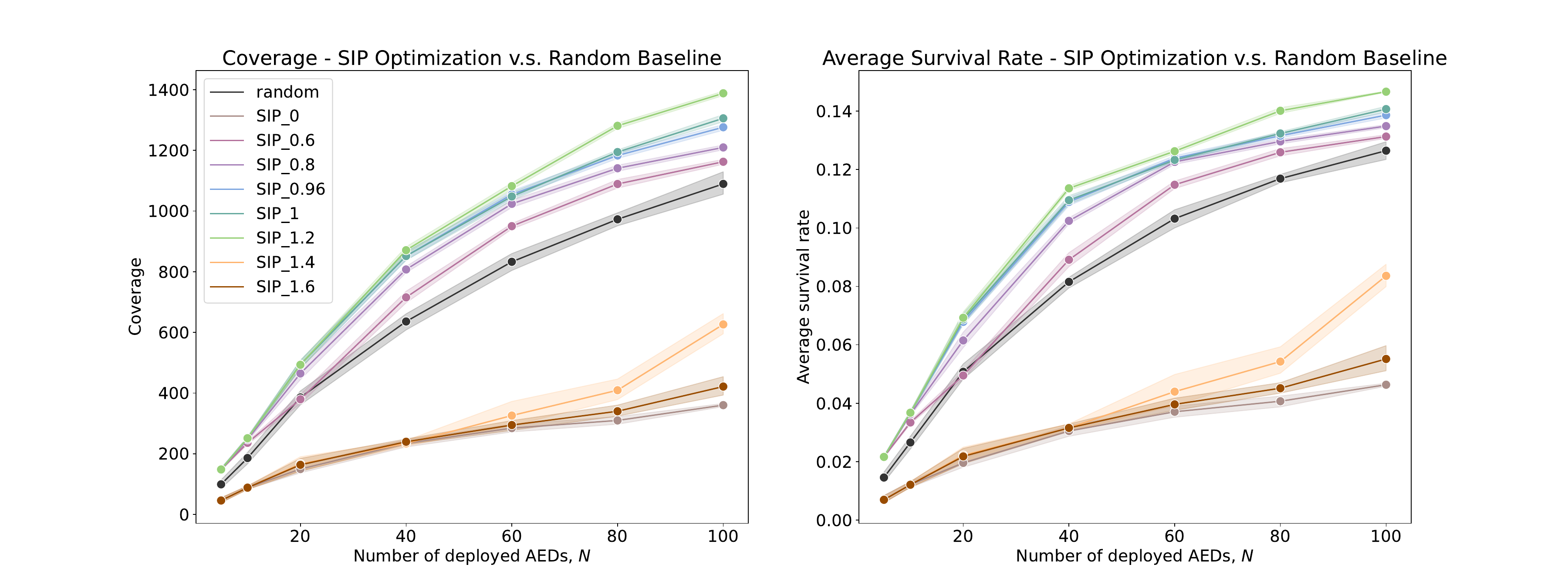}
  \caption{Coverage and average survival rate performance of random baseline and SIP optimization.}
  \label{fig:performance}
\end{figure}

Fig.\ref{fig:performance} compares the coverage and average survival rate of OHCA occurrences between the random baseline and the SIP optimization with varying minimum spacings $D_{min}$.
The detailed incremental percentages are provided in Table \ref{tab:coverage} and \ref{tab:survival_rate}.
The results and observations for both coverage and average survival rate are consistent, indicating a strong correlation between these two metrics: the more OHCA occurrences covered by deployed AEDs, the higher the corresponding average survival rate.
Compared to the random baseline, the SIP solutions, particularly within an appropriate range of $D_{min}$, achieve superior coverage and survival rate performance across all deployment scales $N$ on the real-world OHCA occurrence testing dataset.

For instance, when $D_{min}$ = 1.2 km, the performance of SIP is optimal, outperforming the random deployment baseline by at least 27\% (covering 1388 out of 1504 cases) at $N = 100$, and by up to 49\% (covering 148 cases) at $N = 5$.
Regarding the survival rate, the performance is optimal at $D_{min}$ = 1.2 km, with a minimal increment of 16\%, and up to 48\% at $N = 5$.
In addition to improved performance, the results of SIP models also demonstrate greater robustness, with a lower standard derivation across nearly all cases.
Therefore, the generality of our model's superiority is validated.

\begin{table}[htbp]
\renewcommand\arraystretch{1.3}
    \caption{Coverage performance of SIP optimization under varying minimum spacing.}
    \hspace{-0.5em}
    \begin{tabular}{p{0.4cm}<{\centering} p{1.8cm}<{\centering}
    p{1.7cm}<{\centering} p{0.84cm}<{\centering}
    p{1.7cm}<{\centering} p{0.84cm}<{\centering} 
    p{1.7cm}<{\centering} p{0.84cm}<{\centering} 
    p{1.7cm}<{\centering} p{0.84cm}<{\centering}}
    \toprule
    $D_{min}$ & random & 
    \multicolumn{2}{c}{0} & \multicolumn{2}{c}{0.6} & \multicolumn{2}{c}{0.8} & \multicolumn{2}{c}{0.96} \\
    $N$ & cov. & cov. & \%inc. & cov. & \%inc. & cov. & \%inc. & cov. & \%inc. \\
    \toprule
    5  & 99.1$\pm$21.9  & 45.8$\pm$13.3 & -54\% & 148.1$\pm$5.8 & 49\% & 148.1$\pm$5.8 & 49\% & 148.1$\pm$5.8 & 49\% \\
    10 & 185.8$\pm$31.3 & 88.1$\pm$9.7\ \ \ \  & -53\% & 236.0$\pm$7.9 & 27\% & 247.0$\pm$9.5 & 33\% & 251.1$\pm$7.8 & 35\% \\
    20 & 385.5$\pm$36.7 & 149.4$\pm$21.4 & -61\% & 379.5$\pm$24.2 & -2\% & 464.6$\pm$31.4 & 21\% & 493.3$\pm$21.2 & 28\% \\
    40 & 635.9$\pm$41.4 & 234.3$\pm$19.3 & -63\% & 715.5$\pm$30.7 & 13\% & 807.4$\pm$16.7 & 27\% & 851.8$\pm$22.4 & 34\% \\
    60 & 832.7$\pm$44.1 & 284.4$\pm$20.0 & -66\% & 950.2$\pm$14.4 & 14\% & 1023.8$\pm$19.2 & 23\% & 1054.7$\pm$20.2 & 27\% \\
    80 & 972.8$\pm$35.2 & 309.7$\pm$19.4 & -68\% & 1089.2$\pm$24.4 & 12\% & 1141.1$\pm$16.4 & 17\% & 1183.1$\pm$16.5 & 22\% \\
    100& 1089.5$\pm$60.5 & 360.2$\pm$8.9 & -67\% & 1162.5$\pm$12.7 & 7\% & 1209.7$\pm$12.6 & 11\% & 1276.3$\pm$14.8 & 17\% \\  
    \toprule
    $D_{min}$ & random & 
    \multicolumn{2}{c}{1.0} & \multicolumn{2}{c}{1.2} & \multicolumn{2}{c}{1.4} & \multicolumn{2}{c}{1.6} \\
    \toprule
    5  & 99.1$\pm$21.9 & 148.1$\pm$5.8 & 49\% & 148.1$\pm$5.8 & 49\% & 45.8$\pm$13.3 & -54\% & 45.8$\pm$13.3 & -54\% \\  
    10 & 185.8$\pm$31.3 & 251.1$\pm$7.8 & 35\% & 251.1$\pm$7.8 & 35\% & 88.1$\pm$9.7\ \ \ \  & -53\% & 88.1$\pm$9.7\ \ \ \  & -53\% \\
    20 & 385.5$\pm$36.7 & \textbf{494.4$\pm$21.6} & \textbf{28\%} & 493.2$\pm$22.2 & 28\% & 166.4$\pm$35.6 & -57\% & 163.4$\pm$33.5 & -58\% \\ 
    40 & 635.9$\pm$41.4 & 851.8$\pm$24.2 & 34\% & \textbf{870.9$\pm$18.1} & \textbf{37\%} & 233.9$\pm$15.5 & -63\% & 239.4$\pm$14.8 & -62\% \\
    60 & 832.7$\pm$44.1 & 1048.1$\pm$21.3 & 26\% & \textbf{1082.2$\pm$17.7} & \textbf{30\%} & 325.9$\pm$66.5 & -61\% & 294.5$\pm$26.6 & -65\% \\
    80 & 972.8$\pm$35.2 & 1194.5$\pm$10.8 & 23\% & \textbf{1281.0$\pm$13.9} & \textbf{32\%} & 409.3$\pm$54.6 & -58\% & 340.0$\pm$29.7 & -65\% \\
    100& 1089.5$\pm$60.5 & 1305.7$\pm$19.6 & 20\% & \textbf{1388.2$\pm$11.4} & \textbf{27\%} & 626.3$\pm$51.4 & -43\% & 421.4$\pm$49.6 & -61\% \\  
    \bottomrule
    \multicolumn{10}{l}{note: $D_{min}$: minimum spacing between two AEDs, N: number of total AEDs in the SIP model,}\\
    \multicolumn{10}{l}{\quad \quad \  cov.: coverage, random: random baseline, \%inc.: increase rate compares to random baseline.}
    
    \end{tabular}%
  \label{tab:coverage}%
\end{table}%



\subsubsection{Sensitive Analysis: Deployment Scale $N$ and Minimum Spacing $D_{min}$}

Fig.\ref{fig:performance} demonstrates the intuitive observation that deploying more AEDs leads to higher coverage of OHCA occurrences.
However, the curve increases slower as $N$ increases, indicating a declining marginal effect, suggesting that the deployment of AEDs is approaching saturation when they cover most of the OHCA occurrences (in total 1504).
As observed in Fig.\ref{fig:visualization}, the SIP deployment solution with $N = 100$ has already covered most of the high-risk areas in Virginia Beach.
To reach 100\% coverage, we would need to deploy over 100 additional AEDs throughout the city, with significantly lower efficiency (less than 1 case per AED) compared to the initial $N=100$ AEDs (on average 138 cases per AED).

\begin{figure}[ht!]
  \centering
  \includegraphics[width=\textwidth]{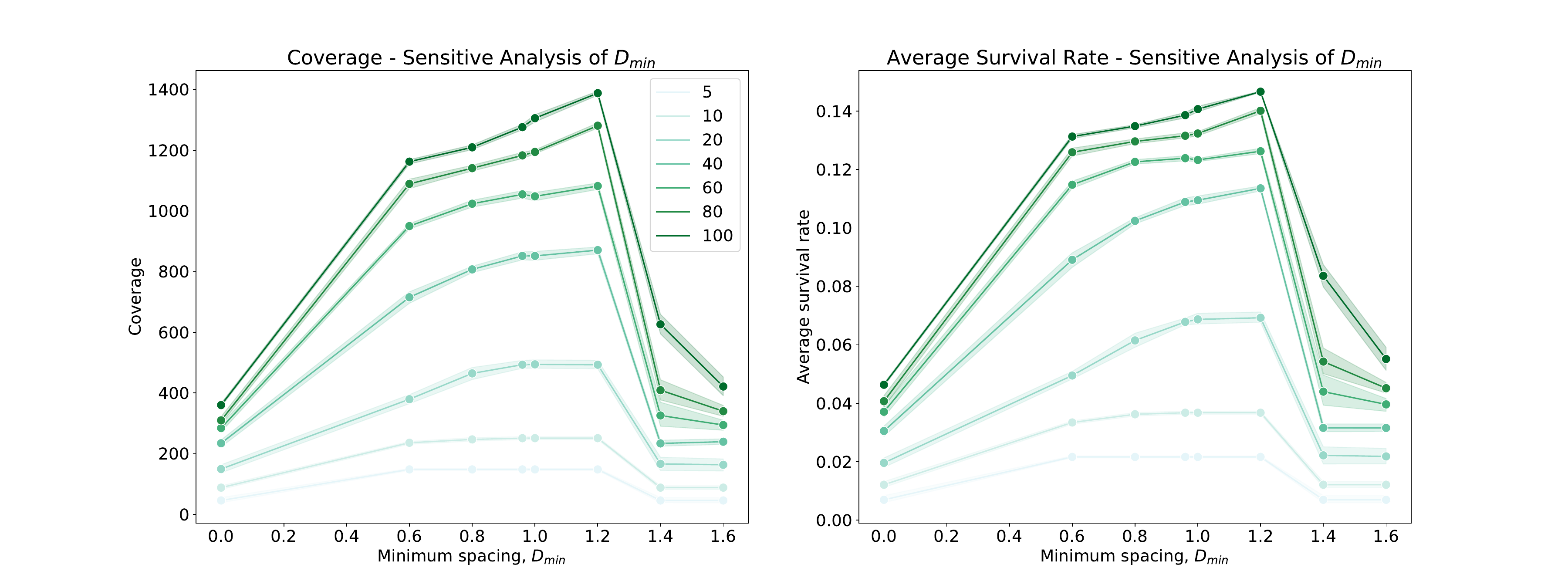}
  \caption{Sensitive analysis of the SIP optimization performance under varying $N$ and $D_{min}$.}
  \label{fig:sensitive_analysis}
\end{figure}

\begin{table}[htbp]
\renewcommand\arraystretch{1.3}
    \caption{Survival rate performance of SIP optimization under varying minimum spacing.}
    \hspace{-1.5em}
    \begin{tabular}{p{0.35cm}<{\centering} p{1.95cm}<{\centering}
    p{1.95cm}<{\centering} p{0.84cm}<{\centering}
    p{1.95cm}<{\centering} p{0.65cm}<{\centering} 
    p{1.95cm}<{\centering} p{0.84cm}<{\centering} 
    p{1.95cm}<{\centering} p{0.84cm}<{\centering}}
    \toprule
    $D_{min}$ & random & 
    \multicolumn{2}{c}{0} & \multicolumn{2}{c}{0.6} & \multicolumn{2}{c}{0.8} & \multicolumn{2}{c}{0.96} \\
    $N$ & sur. & sur. & \%inc. & sur. & \%inc. & sur. & \%inc. & sur. & \%inc. \\
    \toprule
    5  & 0.01$\pm$2.9e-3 & 0.01$\pm$2.0e-3 & -52\% & 0.02$\pm$0.5e-3 & 48\% & 0.02$\pm$0.5e-3 & 48\% & 0.02$\pm$0.5e-3 & 48\% \\
    10 & 0.03$\pm$3.4e-3 & 0.01$\pm$1.6e-3 & -54\% & 0.03$\pm$0.9e-3 & 26\% & 0.04$\pm$1.0e-3 & 36\% & 0.04$\pm$0.9e-3 & 38\% \\
    20 & 0.05$\pm$3.9e-3 & 0.02$\pm$2.6e-3 & -61\% & 0.05$\pm$2.2e-3 & -2\% & 0.06$\pm$4.1e-3 & 21\% & 0.07$\pm$1.8e-3 & 34\% \\
    40 & 0.08$\pm$2.9e-3 & 0.03$\pm$2.8e-3 & -63\% & 0.09$\pm$3.9e-3 & 9\% & 0.1$\pm$1.7e-3 & 26\% & 0.11$\pm$2.2e-3 & 34\% \\
    60 & 0.1$\pm$4.9e-3 & 0.04$\pm$3.0e-3 & -64\% & 0.11$\pm$2.0e-3 & 11\% & 0.12$\pm$1.3e-3 & 19\% & 0.12$\pm$1.4e-3 & 20\% \\
    80 & 0.12$\pm$2.3e-3 & 0.04$\pm$3.0e-3 & -65\% & 0.13$\pm$2.1e-3 & 8\% & 0.13$\pm$1.3e-3 & 11\% & 0.13$\pm$1.4e-3 & 13\% \\
    100& 0.13$\pm$4.9e-3 & 0.05$\pm$0.8e-3 & -63\% & 0.13$\pm$1.1e-3 & 4\% & 0.13$\pm$0.9e-3 & 7\% & 0.14$\pm$1.6e-3 & 10\% \\  
    \toprule
    $D_{min}$ & random & 
    \multicolumn{2}{c}{1.0} & \multicolumn{2}{c}{1.2} & \multicolumn{2}{c}{1.4} & \multicolumn{2}{c}{1.6} \\
    \toprule
    5  & 0.01$\pm$2.9e-3 & 0.02$\pm$0.5e-3 & 48\% & 0.02$\pm$0.5e-3 & 48\% & 0.01$\pm$2.0e-3 & -52\% & 0.01$\pm$2.0e-3 & -52\% \\  
    10 & 0.03$\pm$3.4e-3 & 0.04$\pm$0.9e-3 & 38\% & 0.04$\pm$0.9e-3 & 38\% & 0.01$\pm$1.6e-3 & -54\% & 0.01$\pm$1.6e-3 & -54\% \\
    20 & 0.05$\pm$3.9e-3 & 0.07$\pm$3.0e-3 & 35\% & \textbf{0.07$\pm$2.9e-3} & \textbf{36\%} & 0.02$\pm$4.6e-3 & -56\% & 0.02$\pm$4.3e-3 & -56\% \\ 
    40 & 0.08$\pm$2.9e-3 & 0.11$\pm$2.4e-3 & 34\% & \textbf{0.11$\pm$1.3e-3} & \textbf{39\%} & 0.03$\pm$2.1e-3 & -61\% & 0.03$\pm$2.1e-3 & -61\% \\
    60 & 0.1$\pm$4.9e-3 & 0.12$\pm$0.9e-3 & 19\% & \textbf{0.13$\pm$1.5e-3} & \textbf{22\%} & 0.04$\pm$8.8e-3 & -56\% & 0.04$\pm$3.4e-3 & -62\% \\
    80 & 0.12$\pm$2.3e-3 & 0.13$\pm$1.2e-3 & 13\% & \textbf{0.14$\pm$1.7e-3} & \textbf{20\%} & 0.05$\pm$7.2e-3 & -54\% & 0.05$\pm$2.9e-3 & -61\% \\
    100& 0.13$\pm$4.9e-3 & 0.14$\pm$1.5e-3 & 11\% & \textbf{0.15$\pm$0.3e-3} & \textbf{16\%} & 0.08$\pm$6.1e-3 & -34\% & 0.06$\pm$6.7e-3 & -56\% \\  
    \bottomrule
    \multicolumn{10}{l}{note: $D_{min}$: minimum spacing between two AEDs, N: number of total AEDs in the SIP model,}\\
    \multicolumn{10}{l}{\quad \quad \  sur.: survival rate, random: random baseline, \%inc.: increase rate compared to the random baseline.}
    \end{tabular}%
  \label{tab:survival_rate}%
\end{table}%

Fig.\ref{fig:sensitive_analysis} illustrates the sensitive analysis of performance with respect to the minimum spacing $D_{min}$.
We observe that the SIP model outperforms the random baseline when $D_{min}$ ranges from 0.6 to 1.2 km, with performance significantly declining beyond 1.2 km.
Moreover, the coverage and average survival rate at $D_{min}$ around 1.2 km achieves the minimum standard derivation, indicating that these solutions are the most robust across various candidate sets.
The optimal solution is observed around $D_{min}$ = 1.2 km, which aligns with the $C_{R}$ equal to 0.96 km.
This is the distance that an emergency responder can cover by running from the AED location to the OHCA occurrence site at a speed of 4 m/s within 4 minutes.
The coverage drops as $D_{min}$ decreases because the coverage regions of AED, $A_k$, overlap with each other, leading to inefficient coverage of many OHCA occurrences across the map.
Conversely, as $D_{min}$ increases, AEDs are spaced too far apart to effectively cover the OHCA occurrences in high-risk regions, resulting in a similar decline in coverage efficiency.

\section{Conclusion and Future Work}

Focusing on guaranteeing timely medical interventions for OHCA, this research strives first to make precise and interpretable predictions and then provide accurate and robust AED deployment optimization decisions. To achieve these goals, a novel learn-then-optimize approach is designed, consisting of three key components: a machine learning prediction model, SHAP-based interpretable analytics, and a SHAP-guided integer programming (SIP) model.

In the prediction aspect, considering the difficulty in accessing demographic and historical data, an NN model is trained to utilize only geographic data as inputs and can be generalized to apply in other regions without the need for OHCA historical data. To the best of our knowledge, this study is the first to explore and validate the contribution of geographic data to OHCA occurrences. The predictive performance of our NN model on the test set is verified with $R^2 = 0.752$ and MAE = 5.56. This provides evidence of a strong correlation between geographic features and OHCA occurrences, making it feasible to interpret the NN model with SHAP-based analytics further.

The interpretable SHAP model further quantifies and elaborates on the contribution of each POI or building type to OHCA risk density within a certain region. Based on the average absolute SHAP values, crucial POI and building types have been identified. Moreover, the high positive and negative SHAP values indicate high and low densities of OHCA occurrences, respectively, which align with the reasonable explanations of residential density. These practical insights can help decision-makers in identifying key buildings and regions of AED deployment. Above all, the theoretical findings that the interpretation of geographic data is consistent with the demographic data support the employment of SHAP values in optimizing AED deployment.

An AED deployment optimization integer programming model is established incorporating the computed SHAP-weighted OHCA densities into the objective function. Historical OHCA coverage and average survival rate over historical OHCA incidents are calculated to examine the model optimization effect compared to random baselines. From the results of various numerical experiments on all deployment scales $N$ and an appropriate range of minimum spacing $D_{min}$, our SIP model is verified to present more effective and robust performance over the random baselines, with a 27\%  and 16\% minimum increment in OHCA coverage and survival rate respectively. Sensitive analysis is also conducted to derive some more practical insights. As for deployment scales $N$, the increasing trend of deployment effect along with the increment of $N$ is illustrated and a critical saturated point implying the marginal benefit has been identified. Correspondingly, the AED deployment performance reaches the peak when the minimum spacing $D_{min}$ equals 1.2 km, which demonstrates the significance of minimum spacing selection for decisions in reality.

\deleted[id=Liang]{In conclusion, this study develop an interpretable machine learning model for identifying OHCA high-risk areas and an optimization framework based on SHAP values for AED deployment, aiming to improve accessibility and maximize OHCA survival rates.} In conclusion, this learn-then-optimize approach with SHAP-based interpretable analytics provides a stronger basis for predicting high-risk OHCA areas with easily accessible geographic data and optimizing AED deployment effectively and robustly. It thereby enhances OHCA coverage, improves patients' average survival rates, and strengthens urban emergency service. Furthermore, the proposed approach is not only limited to AED deployment but also can be applied to a broader range of scenarios, such as the distribution of other medical resources like vaccines. 

One potential limitation of this method is that the potential heterogeneity across various cities, such as cultural differences, varying lifestyles, and urban planning, may result in a certain degree of shift in the correlation between geographic features and OHCA occurrences. 
Future research can focus on validating the effectiveness of our proposed approach across different cities.
Alternatively, a benchmark can be developed to assess whether the differences between urban types are significant enough to require re-predictions.
These efforts will enhance the robustness of the model and make it applicable to more cities.

\section*{ACKNOWLEDGMENTS}
This research is funded by the Cheng An Smart Communication Technology (Shenzhen) Co., Ltd. (CACB-2023-A011) and the Guangdong Pearl River Plan (2019QN01X890).

\newpage


\bibliographystyle{wsc}

\begin{thebibliography}{}

\bibitem[\protect\citeauthoryear{Berdowski, Berg, Tijssen, and Koster}{Berdowski et~al.}{2010}]{Berdowski2010}
Berdowski, J., R.~A. Berg, J.~G.~P. Tijssen, and R.~W. Koster. 2010.
\newblock ``Global incidence of out-of-hospital cardiac arrest and survival rates: Systematic review of the literature''.
\newblock {\em Resuscitation\/}~81(11):1479--1487~\url{https://doi.org/10.1016/j.resuscitation.2010.08.006}.


\bibitem[\protect\citeauthoryear{Blom, Beesems, Homma, Zijlstra, Hulleman, Van~Hoeijen, Bardai, Tijssen, Tan, and Koster}{Blom et~al.}{2014}]{blom2014improved}
Blom, M.~T., S.~G. Beesems, P.~C. Homma, J.~A. Zijlstra, M.~Hulleman, D.~A. Van~Hoeijen, , ,  {\em et al}. 2014.
\newblock ``Improved survival after out-of-hospital cardiac arrest and use of automated external defibrillators''.
\newblock {\em Circulation\/}~130(21):1868--1875.


\bibitem[\protect\citeauthoryear{Brodsky}{Brodsky}{2018}]{brodsky2018h3}
Brodsky, I. 2018.
\newblock ``H3: Uber’s hexagonal hierarchical spatial index''.
\newblock {\em Available from Uber Engineering website: https://eng. uber. com/h3/[22 June 2019]\/}:30.


\bibitem[\protect\citeauthoryear{Chen, Dong, Li, Tian, Wu, Li, and Lin}{Chen et~al.}{2024a}]{chen2024optimized}
Chen, H., Y.~Dong, H.~Li, S.~Tian, L.~Wu, J.~Li {\em et al}. 2024a.
\newblock ``Optimized green infrastructure planning at the city scale based on an interpretable machine learning model and multi-objective optimization algorithm: A case study of central Beijing, China''.
\newblock {\em Landscape and Urban Planning\/}~252:105191.


\bibitem[\protect\citeauthoryear{Chen, Dong, Li, Tian, Wu, Li, and Lin}{Chen et~al.}{2024b}]{CHEN2024105191}
Chen, H., Y.~Dong, H.~Li, S.~Tian, L.~Wu, J.~Li {\em et al}. 2024b.
\newblock ``Optimized green infrastructure planning at the city scale based on an interpretable machine learning model and multi-objective optimization algorithm: A case study of central Beijing, China''.
\newblock {\em Landscape and Urban Planning\/}~252:105191~\url{https://doi.org/https://doi.org/10.1016/j.landurbplan.2024.105191}.


\bibitem[\protect\citeauthoryear{Custodio and Lejeune}{Custodio and Lejeune}{2022}]{custodio2022spatiotemporal}
Custodio, J.~E. and M.~A. Lejeune. 2022.
\newblock ``Spatiotemporal data set for out-of-hospital cardiac arrests''.
\newblock {\em INFORMS Journal on Computing\/}~34(1):4--10.


\bibitem[\protect\citeauthoryear{Ekanayake, Meddage, and Rathnayake}{Ekanayake et~al.}{2022}]{ekanayake2022novel}
Ekanayake, I., D.~Meddage, and U.~Rathnayake. 2022.
\newblock ``A novel approach to explain the black-box nature of machine learning in compressive strength predictions of concrete using Shapley additive explanations (SHAP)''.
\newblock {\em Case Studies in Construction Materials\/}~16:e01059.


\bibitem[\protect\citeauthoryear{Ekmekcioglu and Koc}{Ekmekcioglu and Koc}{2022}]{ekmekcioglu2022explainable}
Ekmekcioglu, O. and K.~Koc. 2022.
\newblock ``Explainable step-wise binary classification for the susceptibility assessment of geo-hydrological hazards''.
\newblock {\em Catena\/}~216:106379~\url{https://doi.org/10.1016/j.catena.2022.106379}.


\bibitem[\protect\citeauthoryear{Haklay and Weber}{Haklay and Weber}{2008}]{haklay2008openstreetmap}
Haklay, M. and P.~Weber. 2008.
\newblock ``Openstreetmap: User-generated street maps''.
\newblock {\em IEEE Pervasive computing\/}~7(4):12--18.


\bibitem[\protect\citeauthoryear{Hessulf, Bhatt, Engdahl, Lundgren, Omerovic, Rawshani, Helleryd, Dworeck, Friberg, Redfors, Nielsen, Myredal, Frigyesi, Herlitz, and Rawshani}{Hessulf et~al.}{2023}]{hessulf2023predicting}
Hessulf, F., D.~L. Bhatt, J.~Engdahl, P.~Lundgren, E.~Omerovic, A.~Rawshani, , , , , , , ,  {\em et al}. 2023.
\newblock ``Predicting survival and neurological outcome in out-of-hospital cardiac arrest using machine learning: the SCARS model''.
\newblock {\em EBioMedicine\/}~89:104464~\url{https://doi.org/10.1016/j.ebiom.2023.104464}.
\newblock Epub 2023 Feb 9.


\bibitem[\protect\citeauthoryear{Iban and Sekertekin}{Iban and Sekertekin}{2022}]{iban2022machine}
Iban, M.~C. and A.~Sekertekin. 2022.
\newblock ``Machine learning based wildfire susceptibility mapping using remotely sensed fire data and GIS: A case study of Adana and Mersin provinces, Turkey''.
\newblock {\em Ecological Informatics\/}~69:101647~\url{https://doi.org/10.1016/j.ecoinf.2022.101647}.


\bibitem[\protect\citeauthoryear{Lin, Chu, Lee, and Kao}{Lin et~al.}{2023}]{lin2023optimal}
Lin, C.-H., K.-C. Chu, J.-T. Lee, and C.-Y. Kao. 2023.
\newblock ``Optimal deployment of automated external defibrillators in a long and narrow environment''.
\newblock {\em Plos one\/}~18(2):e0264098.


\bibitem[\protect\citeauthoryear{Lundberg and Lee}{Lundberg and Lee}{2017}]{lundberg2017unified}
Lundberg, S.~M. and S.-I. Lee. 2017.
\newblock ``A unified approach to interpreting model predictions''.
\newblock In {\em Proceedings of the 31st International Conference on Neural Information Processing Systems (NIPS)},  4768--4777.
\newblock California, USA.

\bibitem[\protect\citeauthoryear{Nakashima, Ogata, Kiyoshige, Al-Hamdan, Wang, Noguchi, Shields, Al-Araji, McNally, Nishimura, et~al.}{Nakashima et~al.}{2023}]{nakashima2023machine}
Nakashima, T., S.~Ogata, E.~Kiyoshige, M.~Z. Al-Hamdan, Y.~Wang, T.~Noguchi, , , ,  {\em et al}. 2023.
\newblock ``A machine learning model for predicting out-of-hospital cardiac arrest incidence using meteorological, chronological, and geographical data from the United States''.
\newblock {\em medRxiv\/}:2023--05.


\bibitem[\protect\citeauthoryear{Nichol, Thomas, Callaway, Hedges, Powell, Aufderheide, Rea, Lowe, Brown, Dreyer, Davis, Idris, Stiell, and Investigators}{Nichol et~al.}{2008}]{Nichol2008}
Nichol, G., E.~Thomas, C.~W. Callaway, J.~Hedges, J.~L. Powell, T.~P. Aufderheide, , , , , , ,  {\em et al}. 2008, Sep.
\newblock ``Regional variation in out-of-hospital cardiac arrest incidence and outcome''.
\newblock {\em JAMA\/}~300(12):1423--1431~\url{https://doi.org/10.1001/jama.300.12.1423}.


\bibitem[\protect\citeauthoryear{Nosratabadi, Mosavi, Keivani, Ardabili, and Aram}{Nosratabadi et~al.}{2019}]{nosratabadi2019state}
Nosratabadi, S., A.~Mosavi, R.~Keivani, S.~Ardabili and F.~Aram. 2019.
\newblock ``State of the art survey of deep learning and machine learning models for smart cities and urban sustainability''.
\newblock In {\em International conference on global research and education},  228--238.
\newblock Springer.

\bibitem[\protect\citeauthoryear{Rea, Eisenberg, Sinibaldi, and White}{Rea et~al.}{2004}]{rea2004incidence}
Rea, T.~D., M.~S. Eisenberg, G.~Sinibaldi, and R.~D. White. 2004.
\newblock ``Incidence of EMS-treated out-of-hospital cardiac arrest in the United States''.
\newblock {\em Resuscitation\/}~63(1):17--24.


\bibitem[\protect\citeauthoryear{Ribeiro, Singh, and Guestrin}{Ribeiro et~al.}{2016}]{Ribeiro2016}
Ribeiro, M.~T., S.~Singh, and C.~Guestrin. 2016.
\newblock ``Why should I trust you? Explaining the predictions of any classifier''.
\newblock In {\em Proceedings of the 22nd ACM SIGKDD International Conference on Knowledge Discovery and Data Mining},  1135--1144~\url{https://doi.org/10.1145/2939672.2939778}.

\bibitem[\protect\citeauthoryear{Sasson, Rogers, Dahl, and Kellermann}{Sasson et~al.}{2010}]{sasson2010predictors}
Sasson, C., M.~A. Rogers, J.~Dahl, and A.~L. Kellermann. 2010.
\newblock ``Predictors of survival from out-of-hospital cardiac arrest: a systematic review and meta-analysis''.
\newblock {\em Circ Cardiovasc Qual Outcomes\/}~3(1):63--81~\url{https://doi.org/10.1161/CIRCOUTCOMES.109.889576}.
\newblock Epub 2009 Nov 10.


\bibitem[\protect\citeauthoryear{Sundararajan and Najmi}{Sundararajan and Najmi}{2020}]{sundararajan2020many}
Sundararajan, M. and A.~Najmi. 2020.
\newblock ``The many Shapley values for model explanation''.
\newblock In {\em International conference on machine learning},  9269--9278.
\newblock PMLR.

\bibitem[\protect\citeauthoryear{Yan, Gan, Jiang, Wang, Chen, Luo, Zong, Chen, and Lv}{Yan et~al.}{2020}]{Yan2020}
Yan, S., Y.~Gan, N.~Jiang, R.~Wang, Y.~Chen, Z.~Luo, ,  {\em et al}. 2020, Feb.
\newblock ``The global survival rate among adult out-of-hospital cardiac arrest patients who received cardiopulmonary resuscitation: a systematic review and meta-analysis''.
\newblock {\em Critical Care\/}~24(1):61~\url{https://doi.org/10.1186/s13054-020-2773-2}.


\bibitem[\protect\citeauthoryear{Zheng, Lv, Zheng, Zhang, Tan, Ma, et~al.}{Zheng et~al.}{2023}]{Zheng2023}
Zheng, J., C.~Lv, W.~Zheng, G.~Zhang, H.~Tan, Y.~Ma {\em et al}. 2023.
\newblock ``Incidence, process of care, and outcomes of out-of-hospital cardiac arrest in China: a prospective study of the BASIC-OHCA registry''.
\newblock {\em The Lancet Public Health\/}~8:e923--e932~\url{https://doi.org/10.1016/S2468-2667(23)00139-2}.


\end{thebibliography}

\newpage

\appendix

\section{Appendix: List of POIs and Buildings from OpenStreetMap in Virginia Beach}
\label{appa}
\noindent \textbf{List of POIs: }

Place of worship, grave yard, post office, childcare, courthouse, fire station, library, police, public building, cinema, bar, restaurant, fountain, fast food, cafe, ice cream, dentist, recycling, dojo, pharmacy, atm, clock, parking entrance, bicycle parking, car rental, pub, veterinary, post box, fuel, clinic, bench, bank, parking, social facility, marketplace, sanitary dump station, telephone, nightclub, drinking water, shower, bicycle rental, charging station, loading dock, theatre, community centre, car wash, bicycle repair station, compressed air, letter box, bbq, doctors, planetarium, training, animal boarding, internet cafe, prep school, gambling, driving school, events venue, waste basket, parking space, waste disposal, weighbridge, public bookcase, vending machine, stage, ranger station, animal shelter, exhibition centre, arts centre, kindergarten, bus station, townhall, prison, studio, payment centre.\\

\noindent\textbf{List of Buildings: }

House, bunker, office, commercial, shelter, residential, public, roof, university, dormitory, chapel, greenhouse, apartments, garage, shed, retail, static caravan, church, terrace, service, school, hospital, industrial, pavilion, commercial, stadium, hotel, cabin, toilets, college, warehouse, sports centre, detached, boathouse, barn, riding hall, construction, ship, ruins.

\newpage
\section{Appendix: SHAP Distribution of POIs and Buildings to OHCA Occurrence in each grid in Virginia Beach}
\label{appb}

\begin{figure}[ht!]
  \centering
  \subfigure[Top 38]{
  \includegraphics[width=0.45\textwidth]{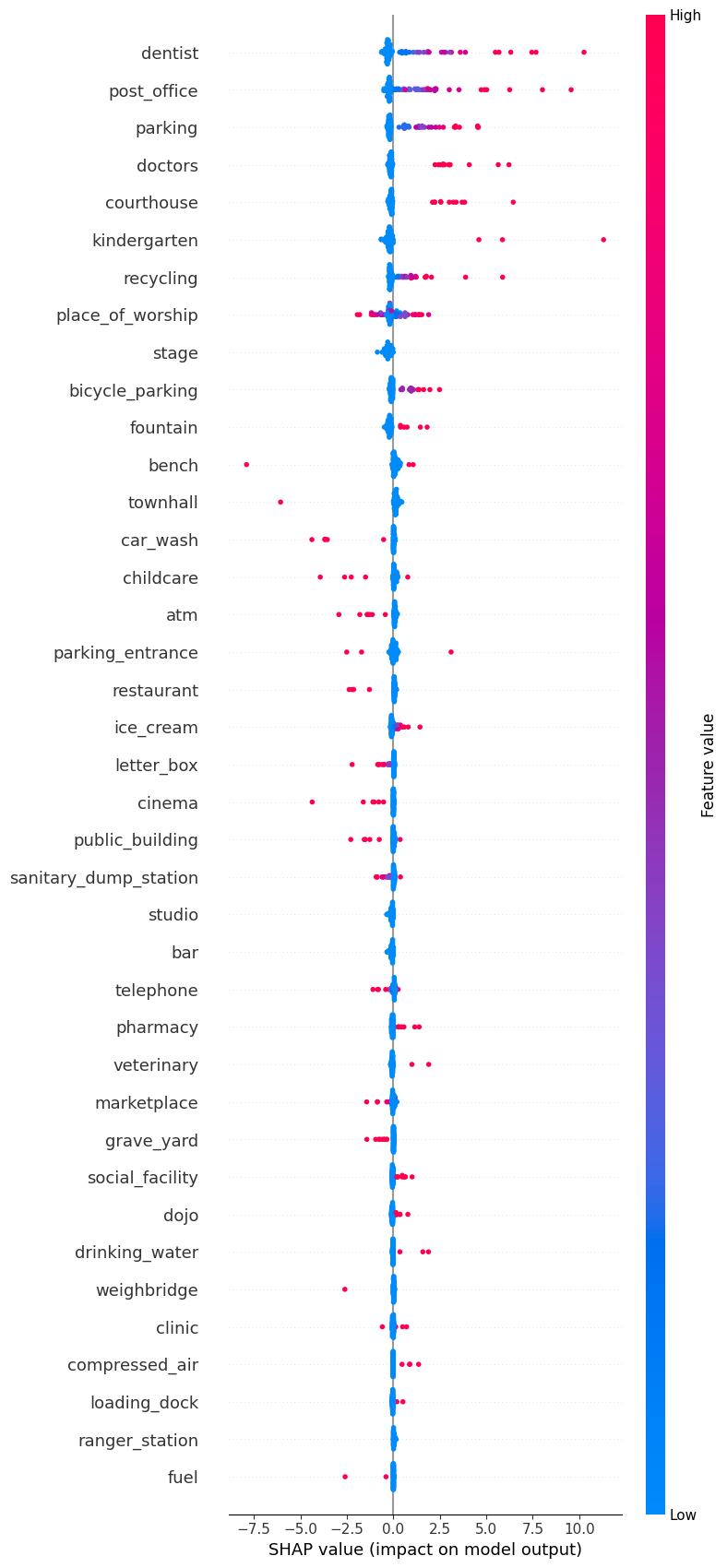}
  \label{fig:shap_app1}
  }
  \subfigure[Last 38]{
  \includegraphics[width=0.45\textwidth]{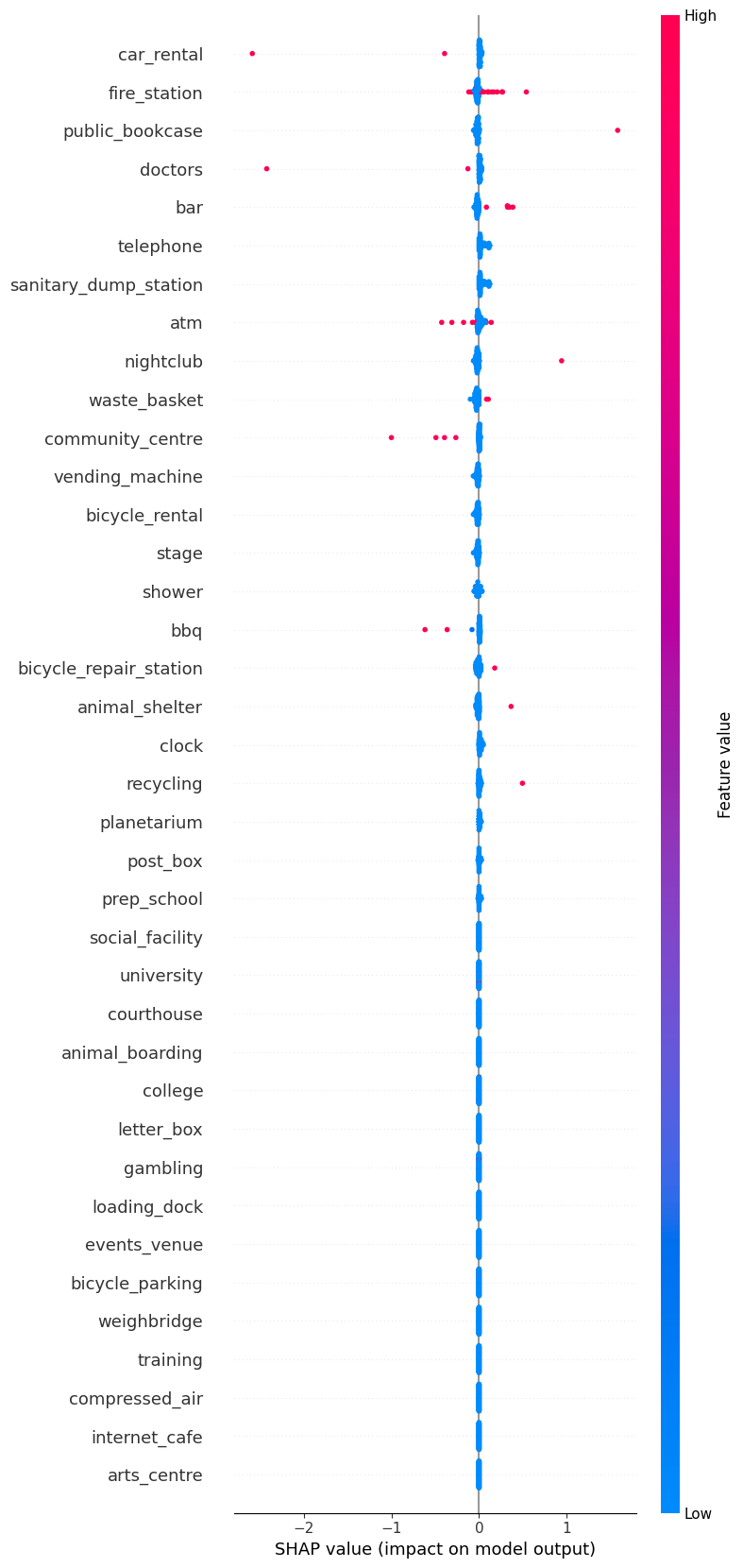}
  \label{fig:shap_app2}
  }

  \caption{SHAP values of all the POIs in the NN model.}
  \label{fig:shap_app}
\end{figure}

\begin{figure}[ht!]
  \centering
  \includegraphics[width=0.5\textwidth]{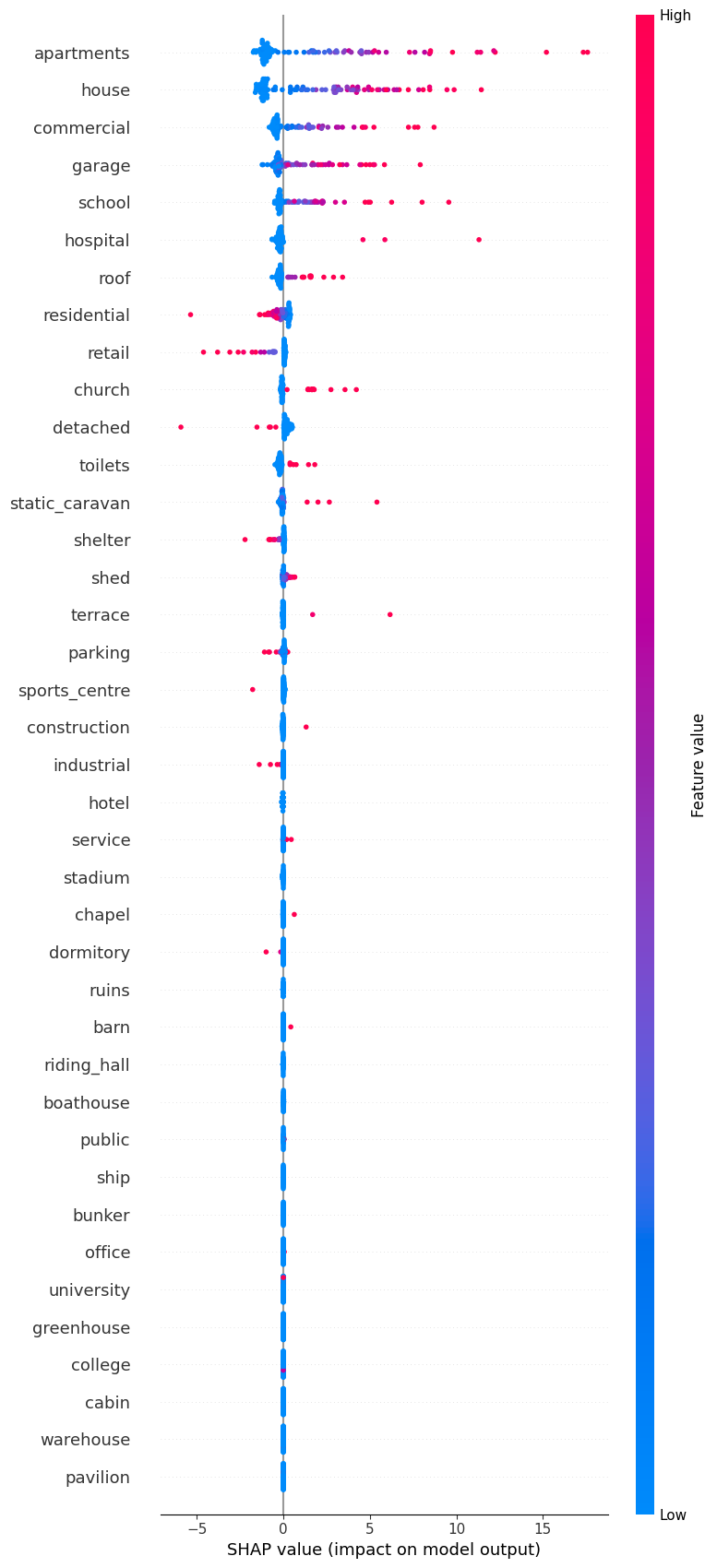}
  \caption{SHAP values of all the buildings in the NN model.}
  \label{fig:shap_app3}
\end{figure}

\end{document}